\documentclass{assets/conf}
\usepackage[inline]{enumitem}
\usepackage{dirtytalk}
\usepackage{subcaption}
\usepackage{lmodern}
\usepackage{amsmath}
\usepackage{multirow}
\usepackage{graphicx}
\usepackage{algpseudocode}
\usepackage{algorithm}
\usepackage{amssymb}
\usepackage{colortbl} 
\usepackage{fancyhdr}
\usepackage{lipsum}
\usepackage{hyperref}
\usepackage{listings}
\usepackage{multicol}
\usepackage{wrapfig}
\usepackage{array}
\usepackage{soul}
\usepackage[dvipsnames]{xcolor}
\usepackage{titlesec}
\titleformat{\section}{\large\bfseries}{\thesection}{1em}{}

\title{Towards Resilient SDA: Graph Theory and Cooperative Control in Distributed Network Architectures}
\author{
{Nesrine Benchoubane}$^{1,2*}$,  
Gunes~Karabulut~Kurt$^{1,2}$ \\
{\textit{\{nesrine.benchoubane, gunes.kurt\}@polymtl.ca}}
}
\institute{\noindent
$^1$Poly-Grames Research Center, Department of Electrical Engineering, Polytechnique Montréal, QC, Canada

$^2$Astrolith, Transdisciplinary Research Unit of Space Resource and Infrastructure Engineering, Polytechnique Montréal, QC, Canada

}
\email{{*CORRESPONDING AUTHOR: Nesrine Benchoubane (e-mail: nesrine.benchoubane@polymtl.ca).}} 

\begin{document}
{\setlength{\parskip}{0pt}
\thispagestyle{plain}
\maketitle 

}

\begin{abstract}
Space Domain Awareness (SDA) involves the detection, tracking, and characterization of space objects through the fusion of data across the space environment. As SDA advances beyond localized or operator-specific capabilities, there is a growing reliance on in-domain space assets for real-time, distributed sensing and decision-making. This paper investigates the potential of on-orbit collaboration by enabling data sharing among heterogeneous satellites as actuators within a single orbital regime. Using graph-theoretic constructs, we define regions of spatial responsibility via Voronoi tessellations and model communication pathways between actuators using Delaunay triangulation. We apply this framework independently to Low Earth Orbit (LEO), Medium Earth Orbit (MEO), Highly Elliptical Orbit (HEO), and Geostationary Orbit (GEO), and analyze each to quantify structural properties relevant to efficient communication, cooperative control, and synchronization for SDA operations with the growth in deployments of space assets.
\end{abstract}


\section*{Nomenclature}
\begin{multicols}{2}
\begin{tabular}{>{}m{1.5cm} m{8cm}}
$C$ & Closeness \\
$D$ & Diameter \\
$H$ & Hopcount \\
$R$ & Radius \\
$a$ & Actuator \\
$A$ & Set of actuators \\
$d$ & Degree \\
$\epsilon$ & Eccentricity \\
$k$ & Coreness \\
$\mathbb{R}^3$ & Three-dimensional space \\
$\mathcal{V}_a$ & Voronoi cell \\
\end{tabular}
\end{multicols}
\section*{Acronyms/Abbreviations}

\begin{multicols}{2}
\begin{tabular}{@{}>{\raggedright\arraybackslash}p{1.5cm} p{5.5cm}@{}}
DT & Delaunay Tetrahedralization \\
FAA & Federal Aviation Administration \\
GEO & Geostationary Orbit \\
GSSAP & Geosynchronous Space Situational Awareness Program \\
HEO & Highly Elliptical Orbit \\
IDAC & Inter-Agency Space Debris Coordination Committee \\
ISL & Inter-Satellite Links \\
LEO & Low Earth Orbit \\
MEO & Medium Earth Orbit \\
MST & Minimum Spanning Tree \\
\end{tabular}

\columnbreak
\begin{tabular}{@{}>{\raggedright\arraybackslash}p{1.5cm} p{5.5cm}@{}}
MW/MT & Missile Warning and Missile Tracking \\
NEOSSat & Near Earth Object Surveillance Satellite \\
ORS-5 & Operationally Responsive Space-5 \\
RSOs & Resident Space Objects \\
SBSS & Space-Based Surveillance Systems \\
SDA & Space Domain Awareness \\
SSN & Space Surveillance Network \\
TLE & Two-Line Element \\
VD & Voronoi Diagram \\
\end{tabular}
\end{multicols}

\section{Introduction}

\begin{figure}
\centering
\includegraphics[width=0.9\linewidth]{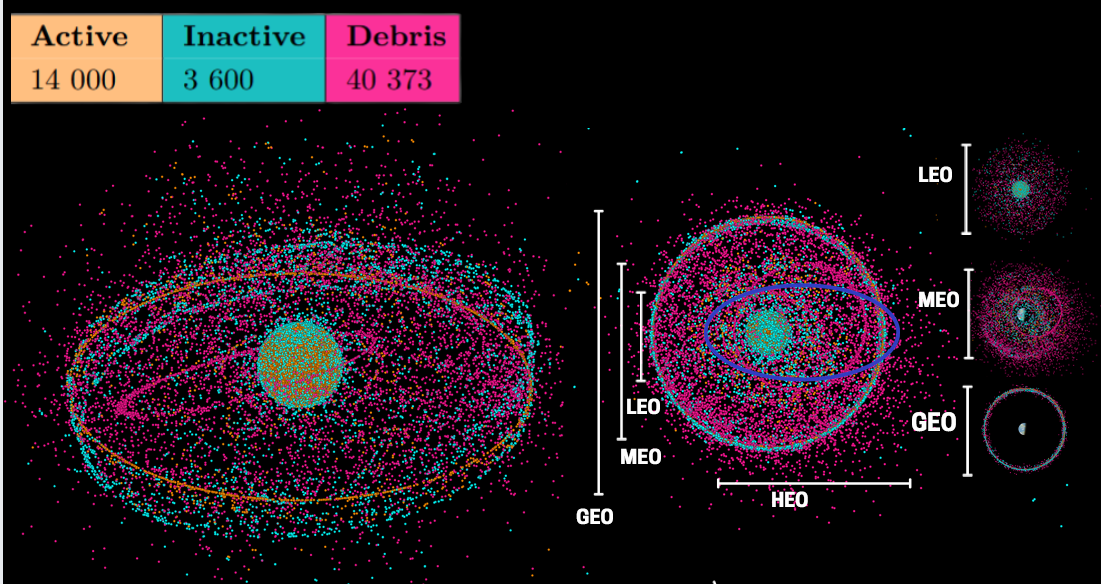}
\captionof{figure}{Distribution of active and inactive RSOs, and debris across orbital regimes (LEO, MEO, GEO, HEO and beyond) based on \cite{ASTRIAGraph, ESA2025SpaceEnvironment}.}
\label{fig:landscape}
\end{figure}

SDA is a critical component of modern space operations, encompassing the detection, tracking, and characterization of RSOs across various orbital regimes. With the accelerated deployment of satellites—driven by the commercial sector, increased accessibility of launch services, and a surge of global actors—the space environment has become more congested, contested, and competitive, as illustrated in Fig. \ref{fig:landscape}. The resulting surge in orbital traffic across LEO, MEO, HEO, GEO, and beyond introduces significant operational challenges for SDA, including collision risks, tracking of non-cooperative or maneuverable objects, and persistent coverage gaps. This growing diversity and volume of space actors, coupled with our own unsustainable practices, have intensified the need for a more cooperative and resilient approach to space stewardship. In recent events, the FAA recently proposed regulations offering U.S. commercial operators multiple disposal options for upper stages, including debris removal contracts, controlled reentry, or transfer to disposal orbits \cite{FAA2023UpperStages}. However, no operational active debris removal systems exist today, despite ongoing development by several companies. The IADC’s most recent status report \cite{IADC2025StatusReport} underscores that adherence to guidelines alone is insufficient, and calls for active remediation strategies to ensure sustainable orbital environments.

SDA operations are further complicated by the increasing prevalence of non-transparent behaviors in space, including unannounced maneuvers, on-orbit interference, and proximity operations that blur the line between benign and hostile intent. These actions, particularly when undertaken by non-cooperative or unregistered entities, increase uncertainty in already congested environments. Traditionally, SDA has relied on ground-based systems—such as the U.S. SSN, which employs radar and electro-optical sensors—face limitations in responsiveness, latency, and full-spectrum observability \cite{Geul2017}. These limitations are \textit{partially} being addressed by space-based systems, including SBSS and ORS-5, which have shown promise in extending SDA capabilities to under-served orbits with improved temporal persistence \cite{Utzmann2014}, but challenges persist in achieving seamless tracking and responsiveness in dynamic orbital regimes \cite{Bloom2022} and continuous monitoring on orbit. 

These latest challenges have spurred significant research leveraging the potential for constellation-level solutions to offer persistent surveillance and improved data fusion across platforms \cite{doi:10.1177/15485129211031673}. Additionally, emerging enabling technologies—including sensor fusion, consensus algorithms, blockchain-secured communication \cite{8898117,benchoubane2025nextgenspacebasedsurveillanceblockchain,9480642, 10908606}, and cooperative control \cite{7536011}—present a path forward. For instance, distributed payloads have been proposed to increase observability and coverage across large-scale constellations \cite{Bertrand2021}, while consensus algorithms are actively being developed to coordinate decisions across heterogeneous on-orbit assets \cite{Klonowski2024} and enhance tracking precision to reduce collision risks\cite{7738363}. 

Building on this momentum, one of the key milestones in the transition to distributed on-orbit architectures has been the publication of latency estimations for distributed SDA in LEO, which demonstrate the feasibility of such systems \cite{gordon2024rolecommunicationsspacedomain}. This breakthrough not only validates the potential of distributed SDA in LEO but also opens the door for new approaches to decentralized tasking and management across other orbital regimes.

\newpage

\noindent \textbf{Motivation and Contributions:} This paper contributes to the shift toward distributed SDA by focusing on intra-orbit collaboration. We propose a graph-theoretic framework to structure communication and task-sharing among satellites. The key contributions of this work are: \begin{enumerate*}[label=(\roman*)]
     \item assessing distributed architectures for SDA across different orbits, with a focus on their synchronization and coordination capabilities through the application of graph theory, and \item developing consensus algorithm tailored to the specific needs of these architectures.
\end{enumerate*}

\section{Distributed On-Orbit SDA}

\begin{figure}[h!]
    \centering
    \includegraphics[width=0.8\linewidth]{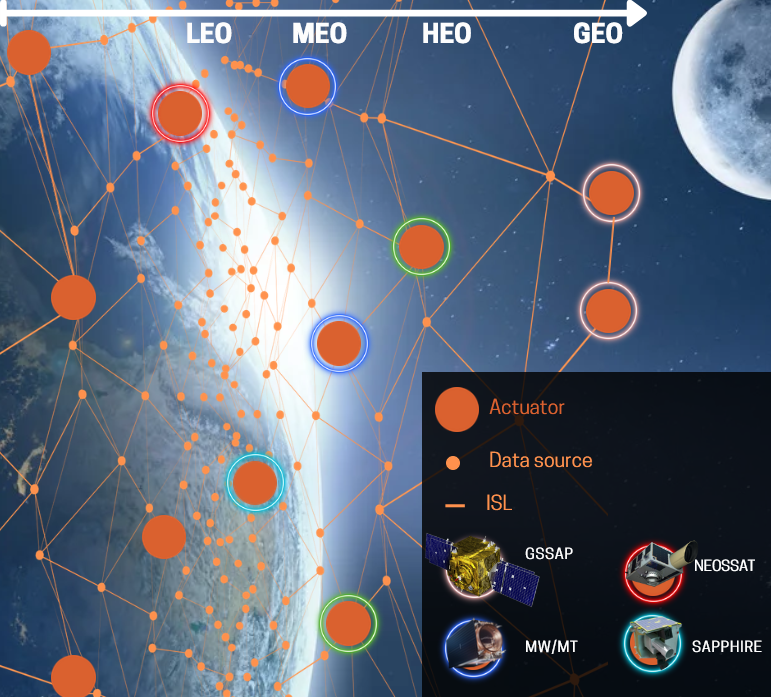}
    \caption{Notional architecture for on-orbit distributed SDA with examples of current and future space surveillance systems such as GSSAP, NEOSSat, MW/MT, and Saphire, assuming roles of actuators with their respective capabilities for enhanced SDA across various orbital regimes.}
    \label{fig:sda-onorbit-inter-orbit}
\end{figure}

In a distributed on-orbit SDA system, as illustrated in Fig. \ref{fig:sda-onorbit-inter-orbit}, the on-orbit nodes can assume two main roles: data source nodes and actuator nodes. Data source nodes are responsible for collecting SDA-related data, such as tracking information, object identification, and space weather parameters. Actuator nodes, in contrast, are equipped with the ability to perform actions based on the data they collect. This on-orbit architecture is also expected to operate in coordination with ground-based assets, together forming a comprehensive, global, and distributed SDA network.

In this paper, we focus primarily on the on orbit architecture with a specific emphasis on the role of actuator nodes. Designing such an architecture requires addressing key operational questions, including: {Which actuator is best suited to communicate with \textit{at a given time?} How can we identify the \textit{closest} actuator to a \textit{specific} data source in real-time? Furthermore, how can we \textit{guide} the efficient coordination of actuators to ensure consistent domain awareness?}


\begin{figure}[h!]
\begin{subfigure}{0.35\linewidth}
\centering
\includegraphics[width=0.8\linewidth]{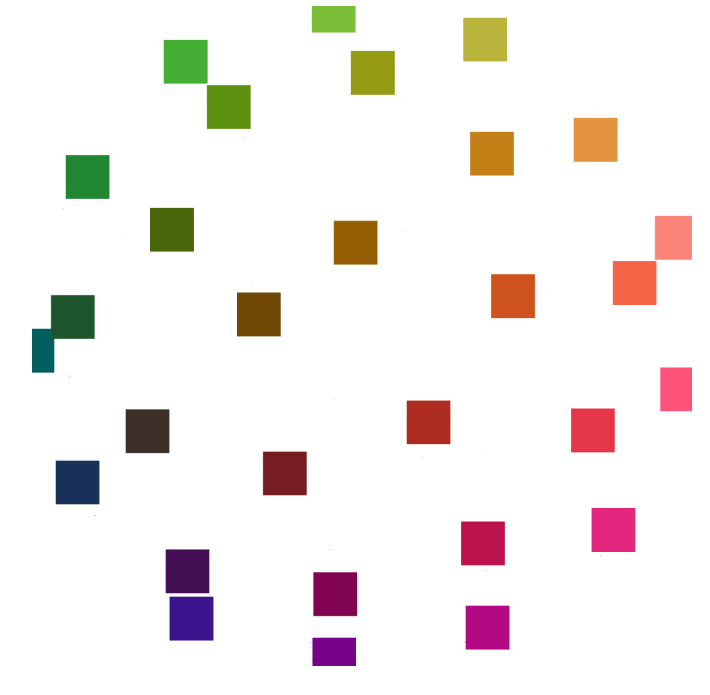}
\caption{}
\label{dvdt-pts}
\end{subfigure}%
\begin{subfigure}{0.35\linewidth}
\centering
\includegraphics[width=0.8\linewidth]{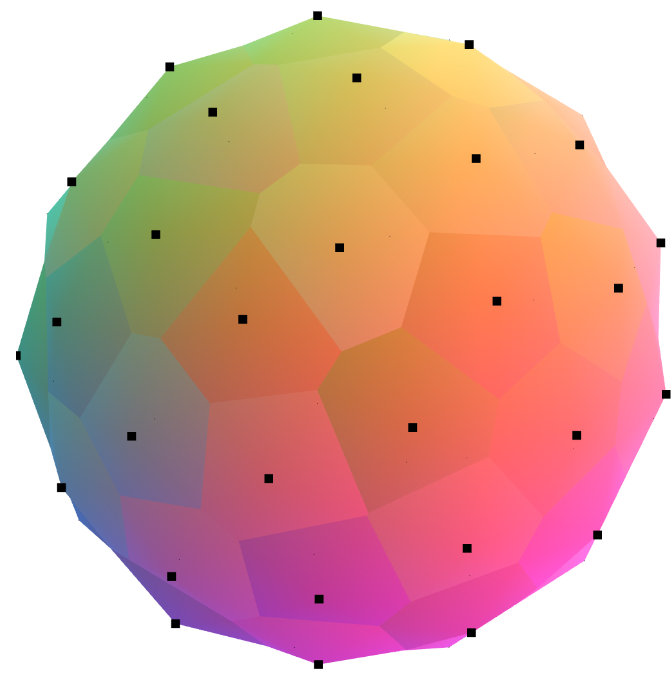}
\caption{}
\label{dv}
\end{subfigure}%
\begin{subfigure}{0.35\linewidth}
\centering
\includegraphics[width=0.8\linewidth]{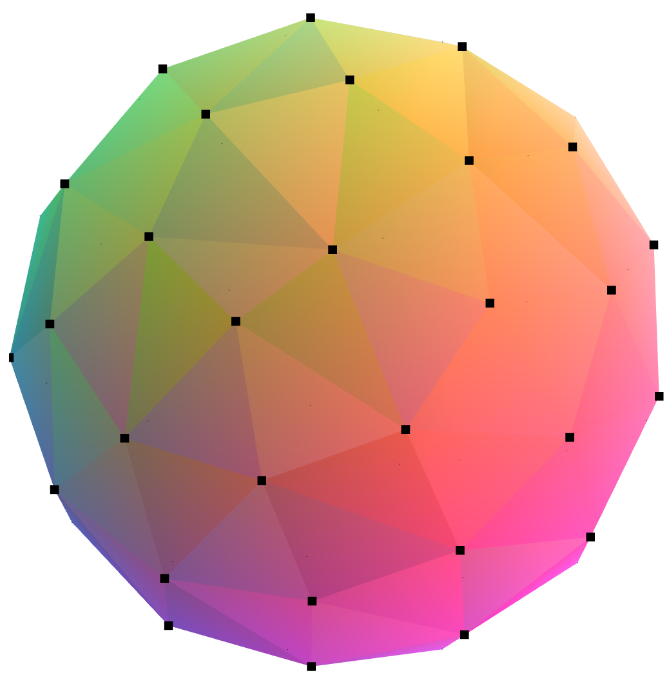}
\caption{}
\label{dt}
\end{subfigure}%
\caption{(a) Distribution of a set of actuators $A=30$ in a sphere. (b) $VD(A)$. (c) $DT(A)$.}
\label{fig:dualdv}
\end{figure}

\subsection{Tessellation for Responsibility Regions}
An effective starting point for managing the operational scope of actuators is spatial partitioning, which divides space into smaller, manageable regions. This approach is well represented by Voronoi tessellation, as illustrated in Fig. \ref{dv}—a technique that partitions space into regions, or Voronoi cells, based on proximity to a predefined set of \say{seeds}. In our context, these seeds represent actuator nodes, and each resulting Voronoi cell defines the region of space that is closest to one actuator relative to all others.

This method naturally lends itself to SDA due to its dynamic, scalable, and geometrically grounded nature. Each actuator is assigned a unique Voronoi cell, enabling it to monitor all entities within its assigned region—whether they are data source satellites, space debris, or other RSOs. This approach has been demonstrated for threat perception regions in \cite{shivshankar2023voronoi}. As actuators move, their associated cells adjust accordingly, maintaining a continuously adaptive map of spatial responsibilities. This ensures that no region is left unmonitored, even as the orbital environment evolves.

An essential feature of Voronoi tessellation is exclusivity: an actuator is associated with exactly one cell at any given time, although this configuration changes dynamically. As the number of actuators increases, the tessellation becomes finer, enhancing resolution and responsiveness—up to certain saturation points, which we investigate in this work.

\subsubsection*{Formal Definition}

Let $A$ be a set of actuators in $\mathbb{R}^d$. The Voronoi cell of an actuator $a \in A$, defined $\mathcal{V}_a$, is the set of all objects $s \in \mathbb{R}^d$ that are closer to $a$ than to any other actuator in $A$; that is,
\begin{equation}
    \mathcal{V}_a = \{ s \in \mathbb{R}^d \mid \|s - a_1\| \leq \|x - a_2\|, \; \forall q \in A \}.
\end{equation}

The union of the Voronoi cells of all generating actuators $a \in A$ form the Voronoi diagram of $A$, defined as ${VD}(A)$. In $\mathbb{R}^3$, ${VD}(A)$, this results in a tessellation with the following key properties:

\begin{itemize}
    \item \textbf{Size:} For $a$ actuators, there are $a$ Voronoi cells in ${VD}(A)$.
    \item \textbf{Voronoi vertices:} They are located at the centers of spheres that touch $4$ actuators in $A$.
    \item \textbf{Voronoi edges:} They are equidistant from $3$ actuators.
    \item \textbf{Voronoi faces:} They are equidistant from $2$ actuators, forming the shared boundary between their regions.
\end{itemize}

\subsection{Mesh-Building for Communication Paths}
While Voronoi tessellation helps define the regions of responsibility for each actuator in space, the next step is to determine communication pathways among them. This is where Delaunay triangulation—the geometric dual of the Voronoi diagram—becomes essential. As illustrated in Fig.~\ref{dt}, Delaunay triangulation connects actuators that share a boundary in the Voronoi diagram. 

The Delaunay structure is proximity-driven: it links actuators that are spatially adjacent by drawing edges between them. These edges indicate potential communication channels, thereby defining a topology that is both adaptive and optimized for locality. 

\subsubsection*{Formal Definition}

Let $A$ be a set of actuators in $\mathbb{R}^d$. The Delaunay triangulation $DT(A)$ is a triangulation of the convex hull of $A$ such that no actuator in $A$ lies inside the circumsphere of any simplex in the triangulation,
\begin{equation}
    \forall \sigma \subset A, \quad \text{circumsphere}(\sigma) \cap A = \sigma.
\end{equation}

That is, each simplex in the Delaunay triangulation has an empty circumsphere—no other actuator lies within the sphere defined by its vertices, which ensures that we do not create edges that would span over a closer actuator. This ensures the locality mentioned before. This also results with following interesting key properties:

\begin{itemize}
    \item \textbf{Euclidean Minimum Spanning Tree:} : The MST of the actuator set $A$ is \textit{always} a subgraph of $DT(A)$. This means that the Delaunay triangulation captures the minimal infrastructure needed for connectivity.
    \item \textbf{Spanner Property}: Delaunay edges serve as a good geometric spanner, offering near-optimal paths between actuators. That is, the path between two actuators in $DT(A)$ is not much longer than the direct Euclidean distance, making it suitable for efficient routing.
    \item \textbf{Closest Pair Property}: The closest pair of actuators in $A$ are guaranteed to be directly connected in the Delaunay triangulation. This ensures that immediate neighbors are always part of the communication mesh.
\end{itemize}






\subsection{Graph Metrics}
\label{sec: graph}


To assess the performance of the networks formed by the dual Voronoi-Delaunay graphs over time, we rely on graph-based metrics that can be categorized into two main groups: \begin{enumerate*}[label=(\roman*)] \item \textit{topology-based metrics}, which are governed by the spatial distribution and relative arrangement of the nodes in the graph; and \item \textit{service-based metrics}, which capture operational characteristics of the nodes as well as the links.\end{enumerate*}

Since the graph's structure is dictated by the communication pathways between actuators based on proximity and spatial arrangement, we propose evaluating the network's performance through topological graph metrics which include distance-based as well as connectivity-based metrics. These metrics help us understand various aspects of the network configuration and its implications, including:
\begin{itemize}
    \item Communication efficiency, by quantifying how quickly information can propagate between actuators.
    \item Synchronization potential, by examining how central or peripheral nodes are in the overall network.
    \item Cooperative control feasibility, by identifying strong substructures or weak links that impact collective behaviors and distributed decision-making.
\end{itemize}

These topological metrics thus offer both diagnostic and design value, helping identify bottlenecks, optimize actuator placement, and assess how structural changes could influence the goal we set for the system. While our current focus is on topology-based metrics, these findings will later serve as a scaffold for service-level evaluations that incorporate the operational characteristics. 

\subsubsection{Distance metrics} These metrics characterize the shortest-path properties of the network and directly influence communication latency and the ability of actuators to synchronize actions.

\paragraph{Hopcount:} For any two actuators $a_i$ and $a_j$, the hop count $H_{a_i \to a_j}$ is defined as the number of edges in the shortest path between them,
\begin{equation}
    H_{a_i \to a_j} = min_{k \in [1, A-1]} (P_{{a_i} \to a_j} (k)),
\end{equation}

Furthermore, the distribution of hop counts for different $k$ values allows us to assess the expected number of hops in common communication scenarios, helping to evaluate how efficient or costly routing might be in sparse vs. dense deployments.

\begin{quote}
\textbf{Rationale:} \textit{Hop count is directly related to latency. By minimizing the average and worst-case hop counts, we can improve the end-to-end delivery between the two most distant actuators.}
\end{quote}

\paragraph{Closeness:} For a given actuator, on average, the closeness measures how close it is to all the rest in the network. For an actuator $a_i$, 
\begin{equation}
    C_i = \Bigg({\sum_{a_j \in {A} \setminus \{a_i\}} H_{a_i \to a_j}} \Bigg)^{-1}.
\end{equation}

\begin{quote}
    \textbf{Rationale:} \textit{High closeness score suggests that an actuator can disseminate information quickly and is likely to play a critical role in time-sensitive mission operations. Actuators with lower scores are more isolated and may experience higher.}
\end{quote}

\paragraph{Eccentricity:} For a given actuator, $\epsilon_i$ reflects the worst-case communication delay from there to any other in the network. For an actuator $a_i$, 
\begin{equation}
    \epsilon_i = \text{max}_{a_j \in {A}} (H_{a_i \to a_j}).
\end{equation}

\begin{quote}
    \textbf{Rationale:} \textit{Actuators with lower eccentricity are more central and better suited for coordinating global tasks. Eccentricity also reveals structural bottlenecks and potential communication vulnerabilities.}
\end{quote}

\paragraph{Diameter:} For the network, $D$ is the maximum of all actuator eccentricities and represents the longest shortest path in the network,
\begin{equation}
    D = \text{max}_{a_i \in {A}} (\epsilon_i).
\end{equation}

\begin{quote}
    \textbf{Rationale}: \textit{A larger diameter suggests that information could take longer to propagate through the network and indicate long chains of actuators.
}
\end{quote}

\paragraph{Radius:} For the network, $R$ is the minimum of all actuator eccentricities and indicates how close the most "central" actuator is to the rest,
\begin{equation}
    R = \text{min}_{a_i \in {A}} (\epsilon_i).
\end{equation}

\begin{quote}
    \textbf{Rationale}: \textit{It helps identify highly accessible actuators in the network and the most central in terms of reachability.}
\end{quote}

\subsubsection{Connection metrics} These metrics capture the direct interaction patterns and structural robustness of the network and consider the cohesive groups that can be formed from within; these would be subgroups among actuators that have strong, direct, and \textit{frequent} ties.

\paragraph{Node Degree: } For a given actuator, $d$ refers to the number of immediate connections it has to other actuators,
\begin{equation}
    d_i = \sum_{j=1}^{A} a_{ij},
\end{equation}
where $a_{ij} = 1$ if a direct link exists between $a_i$ and $a_j$, and 0 otherwise if no links (In our case, this doesn't exist due to the properties of the DT). In addition, the average and variance of the node degree across the graph is evaluated to provide insights into network redundancy and balance.

The degree distribution for a set value $k$ is useful for understanding resilience:
\begin{equation}
    Pr[D=k] = \frac{d_k}{A},
\end{equation}
where $d_k$ is the number of actuaors with degree $k$. 

\begin{quote}
    \textbf{Rationale}: \textit{High-degree actuators facilitate rapid information spread and coordination. The shape of the distribution also reflects network resilience.}
\end{quote}

\paragraph{Coreness ($k$):} For the network, the $k$-core is the maximal subgraph in which every actuator has at least degree $k$. Coreness of an actuator is the highest$k$for which it belongs to a $k$-core.

\begin{quote}
    \textbf{Rationale}: \textit{High coreness implies structural centrality and deep embedding within the network. Such actuators are stable candidates for leading synchronization tasks or initiating cooperative control due to their topological influence.}
\end{quote}

\section{Simulation Framework for Orbital Actuator Distribution}

We retrieve orbital data from Space-Track \cite{spacetrack} in TLE format, which we then propagate into space-time coordinates using the \textsl{satellite-js} library \cite{satellitejs}. These satellites serve as proxies for actuator positions, grouped by their respective orbits, as follows: LEO: 2,564, MEO: 183, GEO: 839, and HEO: 182. Our aim is to simulate the scale of these constellations as well as different orbital regime distributions.

\vspace{0.3cm}
\begin{lstlisting}[frame=single]
0 NEOSSAT  

1 39089U 13009D   18115.60454839 +.00000035 +00000-0 +27923-4 0  9992  

2 39089 098.5303 320.4424 0012356 092.9866 267.2733 14.34421818270178  
\end{lstlisting}

For instance, the TLE entry above describes NEOSSat with three key lines: \begin{enumerate*}[label=(\roman*)] \item Line 1 is the satellite's name. \item Line 2 contains the NORAD ID, classification, epoch, and drag info. \item Line 3 encodes orbital elements: inclination, RAAN, eccentricity, argument of perigee, mean anomaly, and mean motion.\end{enumerate*}

\begin{algorithm}[h!]
\caption{Generation of the VD and DT over time.}\label{alg:cap}
\begin{algorithmic}
\Require An orbital regime, a set of actuators $A$, a percentage of actuators to extract 
\Require A list of TLE elements for all extracted actuators
\Require Time step for periodic updates
\Ensure Updated $VD(A)$, $DT(A)$, and actuators positions over time

\State Initialize Earth projection with \texttt{geoOrthographic()} for visualization
\State Set time for simulation start (e.g., current timestamp)
\State Initialize data storage (CSV for results)
\State Set the duration for simulation to 24 hours

\While{simulation time < 24 hours}
    \State Propagate actuators positions based on TLE data and current timestamp
    \For{each actuator}
        \State Compute actuator's position at current time using \texttt{actuator.propagate()}
        \State Convert position to Earth-centered geodetic coordinates
        \State Store actuator's position and timestamp
    \EndFor
    \State Generate $VD(A)$ diagram based on actuator positions
    \State Generate $DT(A)$ from Voronoi diagram
    \State Save current positions, $VD(A)$, $DT(A)$ to CSV
    \State Update simulation time by adding time step
\EndWhile

\State Finalize CSV export with complete actuator data and topologies.
\end{algorithmic}
\end{algorithm}

\begin{figure}[h!]
\centering
  \begin{subfigure}{0.49\linewidth}
    \centering
    \includegraphics[width=\linewidth]{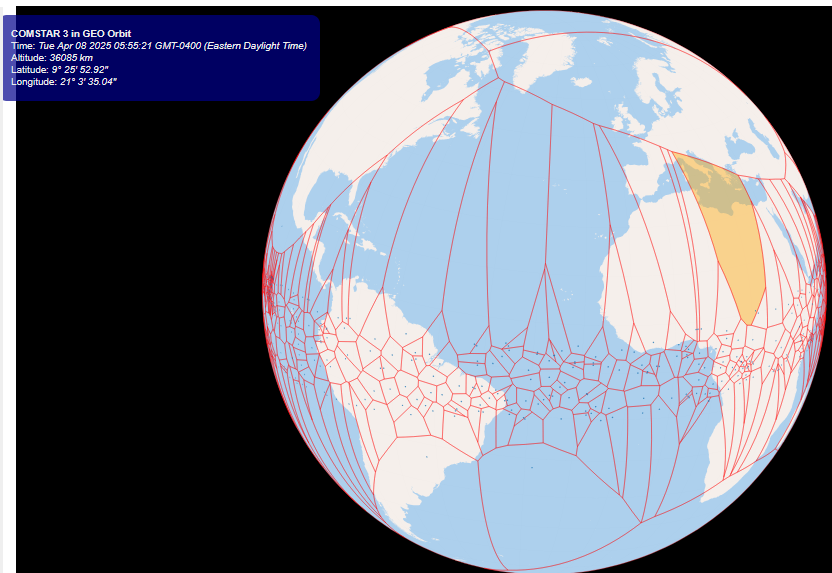}
    \caption{GEO}
  \end{subfigure}%
  \begin{subfigure}{0.5\linewidth}
  \centering
    \includegraphics[width=\linewidth]{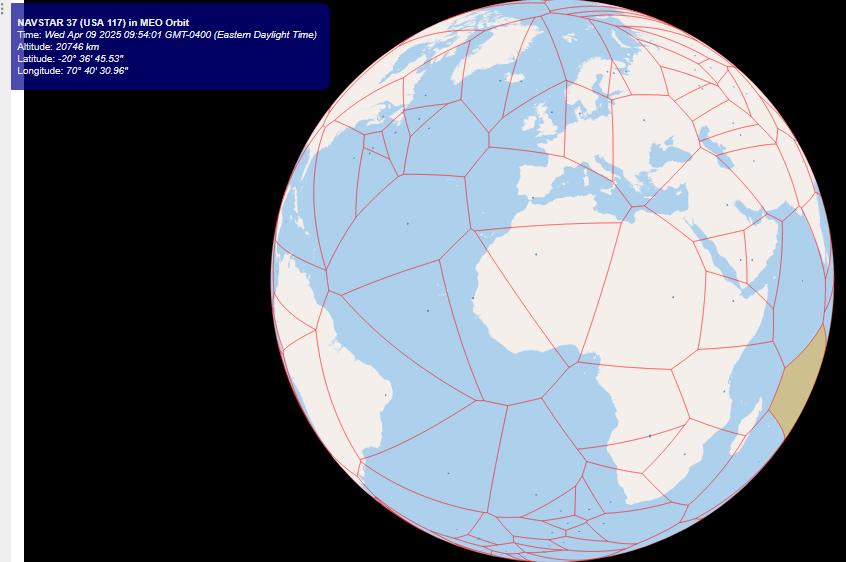}
    \caption{MEO}
  \end{subfigure}
  \caption{Voronoi diagrams illustrating spatial coverage of actuators in GEO and MEO.}
    \label{fig:diagram}
\end{figure}

To compute the Voronoi and Delaunay diagrams, the actuator positions are projected to their respective zenith points on Earth. These projections are then used to generate Voronoi diagrams and corresponding Delaunay triangulations using the \textsl{d3-geo-voronoi} library \cite{d3geovoronoi}. The input for the Voronoi computation consists of arrays of [longitude, latitude] coordinates for each actuator. For context, we overlay the actuator network on 1:10m Cultural Vector data from Natural Earth \cite{naturalearth_admin0}, following the visualization framework described in \cite{vasturiano_voronoi}.

Next, we extract 24-hour intervals of the generated topologies, with a snapshot shown in Fig. \ref{fig:diagram}, to capture the temporal evolution of the network. Algorithm~\ref{alg:cap} describe all above operations where we vary the initial actuator participation from 5\% to 95\% to study how partial deployment influences the network structure and resilience over time. These time-series topologies include altitude, longitude, latitude, satellite name, orbit type, and timestamps, and are derived from Delaunay triangulations. The resulting data is exported and analyzed using the metrics outlined in Section \ref{sec: graph}.

\section{Analysis and Discussion per Orbital Regime}



\begin{table}[ht]
\centering
\caption{Summary of graph metrics across orbits for varying percentages of actuators.}
\label{tab:summ}
\renewcommand{\arraystretch}{1.2}
\begin{tabular}{|c|c|c|c|c|c|c|}
\hline
\textbf{Orbit} & \textbf{\% Actuators} & \textbf{Mean $\epsilon$} & \textbf{Max $\epsilon$} & \textbf{Mean $D$} & \textbf{Mean $R$} & \textbf{Min $R$} \\
\hline\hline
\multirow{5}{*}{LEO} 
& 10 & 4.000 & 4 & 4.000 & 2.000 & 2 \\
& 30 & 3.951 & 4 & 3.951 & 2.000 & 2 \\
& 50 & 3.672 & 4 & 3.672 & 2.000 & 2 \\
& 70 & 3.246 & 4 & 3.246 & 2.000 & 2 \\
& 90 & 2.852 & 4 & 2.852 & 1.869 & 1 \\
\hline\hline
\multirow{4}{*}{MEO} 
& 10 & 3.096 & 4 & 3.096 & 2.000 & 2 \\
& 30 & 4.652 & 6 & 4.652 & 2.888 & 2 \\
& 50 & 5.175 & 6 & 5.175 & 3.045 & 2 \\
& 70 & 5.726 & 7 & 5.726 & 3.242 & 2 \\
\hline\hline
\multirow{4}{*}{HEO} 
& 10 & 2.941 & 4 & 2.941 & 2.000 & 2 \\
& 30 & 4.037 & 5 & 4.037 & 3.000 & 3 \\
& 50 & 4.906 & 5 & 4.906 & 3.350 & 3 \\
& 70 & 5.063 & 6 & 5.063 & 3.998 & 3 \\
\hline\hline
\multirow{4}{*}{GEO} 
& 10 & 5.513 & 6 & 5.513 & 3.000 & 3 \\
& 30 & 6.653 & 8 &  6.653 & 4.000 & 4 \\
& 50 & 7.231 & 8 &  7.231 & 4.560 & 4 \\
& 70 & 8.276 & 9 & 8.276 & 4.997 & 4 \\
\hline
\end{tabular}
\end{table}

We begin by analyzing key distance metrics across the orbital regimes which are summarized in Table~\ref{tab:summ} and represent averaged values over the simulation time frame for each orbital regime and varying actuator percentages. In LEO, the network structure demonstrates improved performance as more actuators are added. The mean graph eccentricity and diameter both decrease steadily from 4.0 at 10\% to approximately 2.85 at 90\%, which that as more actuators are introduced, the graph becomes more compact, leading to shorter communication paths. The radius mean remains constant at 2.0 until it drops slightly to 1.87 at 90\%, and the radius minimum also drops to 1 which overall demonstrate that the LEO networks scale efficiently and benefit from increased actuators density, maintaining low communication latency and strong overall connectivity. However, MEO shows a reverse trend, where network compactness deteriorates as actuators increases. The mean graph eccentricity rises from approximately 3.1 at 10\% to 5.7 at 70\%, and the graph diameter follows a similar pattern; which indicate that the network becomes less centralized and more fragmented at higher actuators density levels. This degradation suggests that MEO constellations are more susceptible to structural inefficiencies under increased node density, due to larger inter-satellite distances. Similar pattern is observed in HEO and GEO, beginning with relatively efficient connectivity at 10\% in their respective orbit which deteriorate as actuator percentage increases.

\begin{figure}[h!]
\centering
\begin{subfigure}{0.31\linewidth}
  \centering
    \includegraphics[width=\linewidth, trim=0 0 0 40, clip]{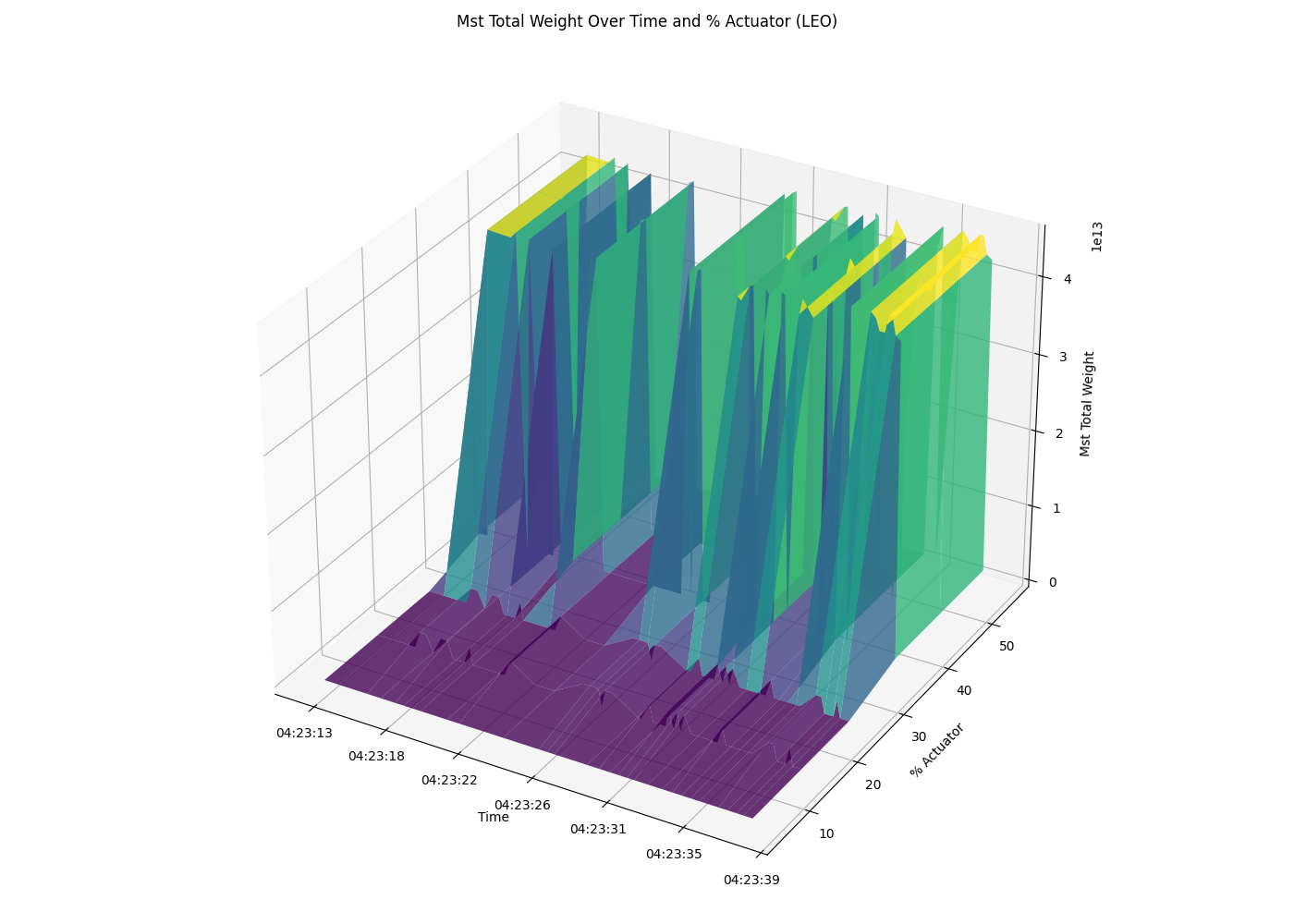}
\caption{}
\end{subfigure}%
\begin{subfigure}{0.31\linewidth}
  \centering
    \includegraphics[width=\linewidth, trim=0 0 0 40, clip]{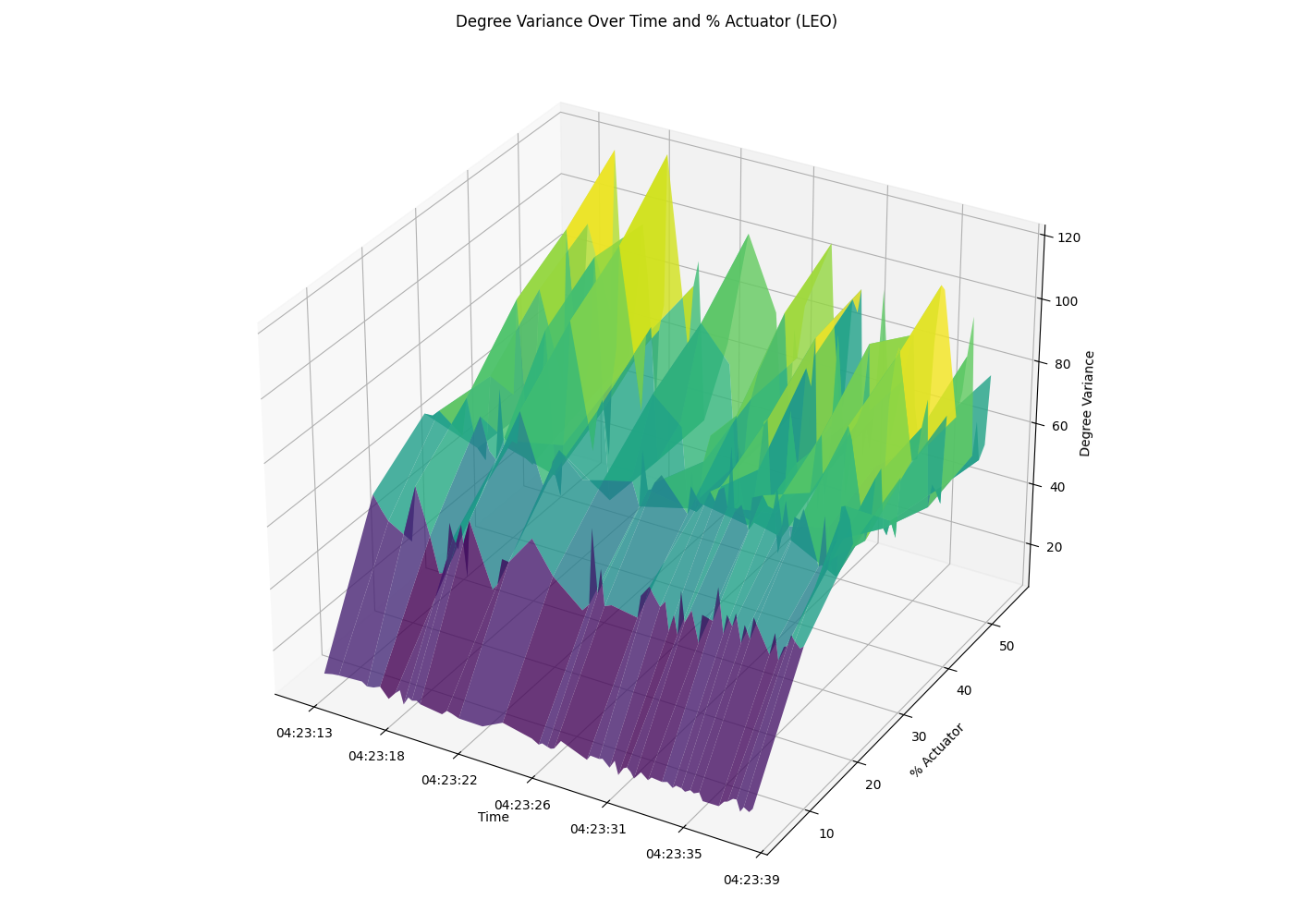}
\caption{}
\end{subfigure}%
\begin{subfigure}{0.31\linewidth}
  \centering
    \includegraphics[width=\linewidth, trim=0 0 0 40, clip]{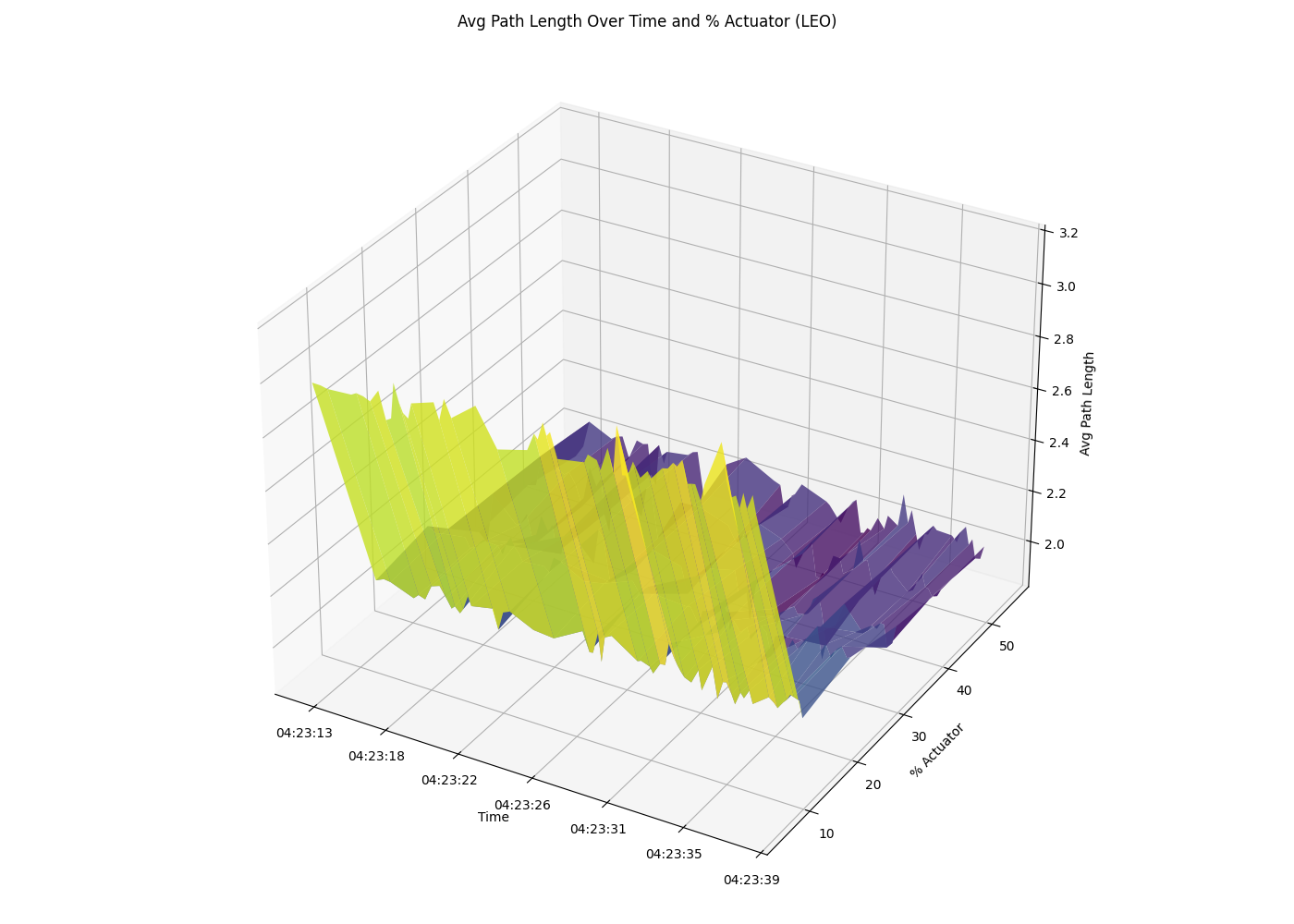}
\caption{}
\end{subfigure}
\caption{Performance metrics for LEO actuator network topology over time. (a) MST, (b) Degree variance, and (c) Average path length.}
\label{fig:leo2}
\end{figure}

\begin{figure}[h!]
\centering
\begin{subfigure}{0.31\linewidth}
  \centering
    \includegraphics[width=\linewidth, trim=0 0 0 40, clip]{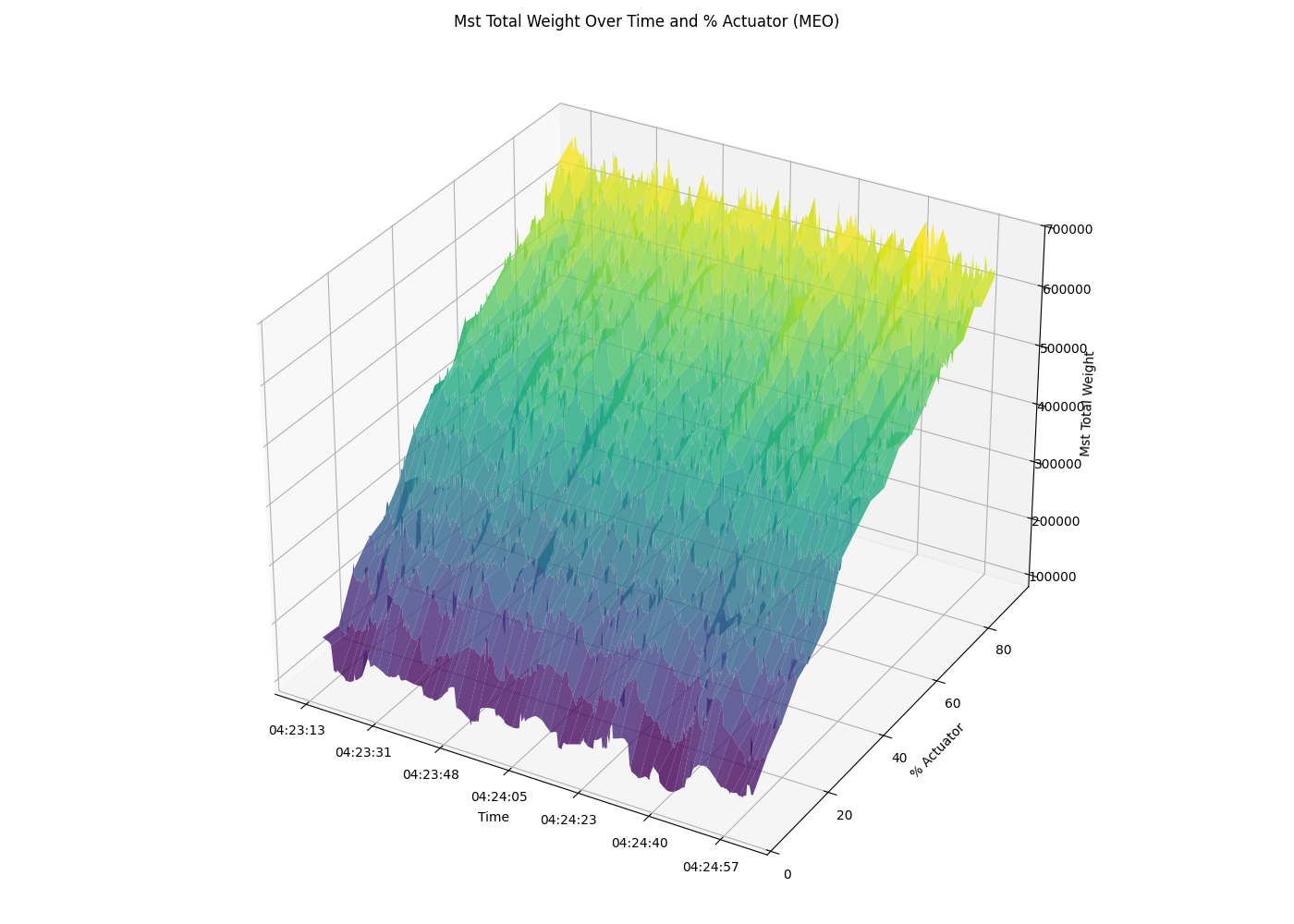}
\caption{}
\end{subfigure}%
\begin{subfigure}{0.31\linewidth}
  \centering
    \includegraphics[width=\linewidth, trim=0 0 0 40, clip]{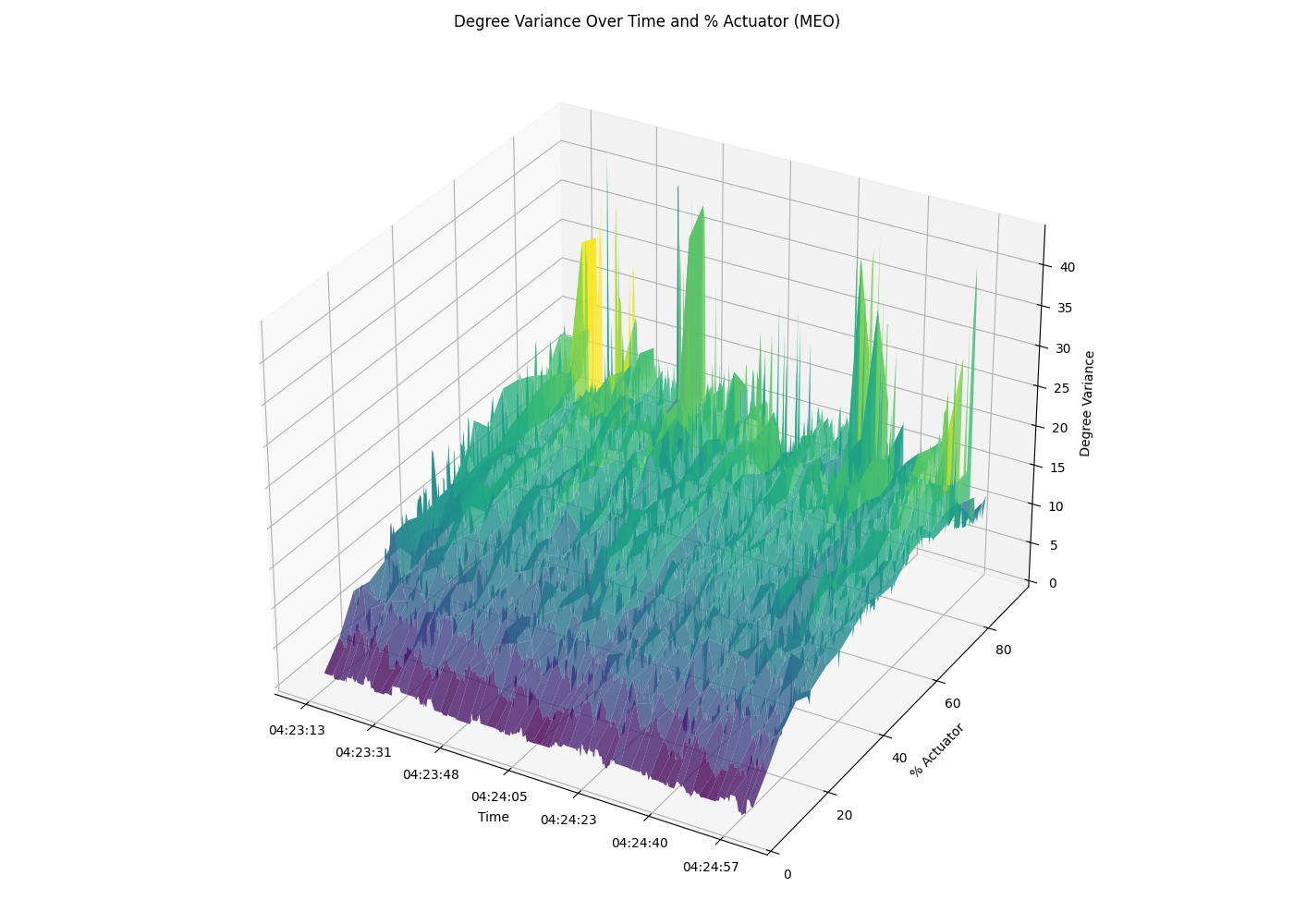}
\caption{}
\end{subfigure}%
\begin{subfigure}{0.31\linewidth}
  \centering
    \includegraphics[width=\linewidth, trim=0 0 0 40, clip]{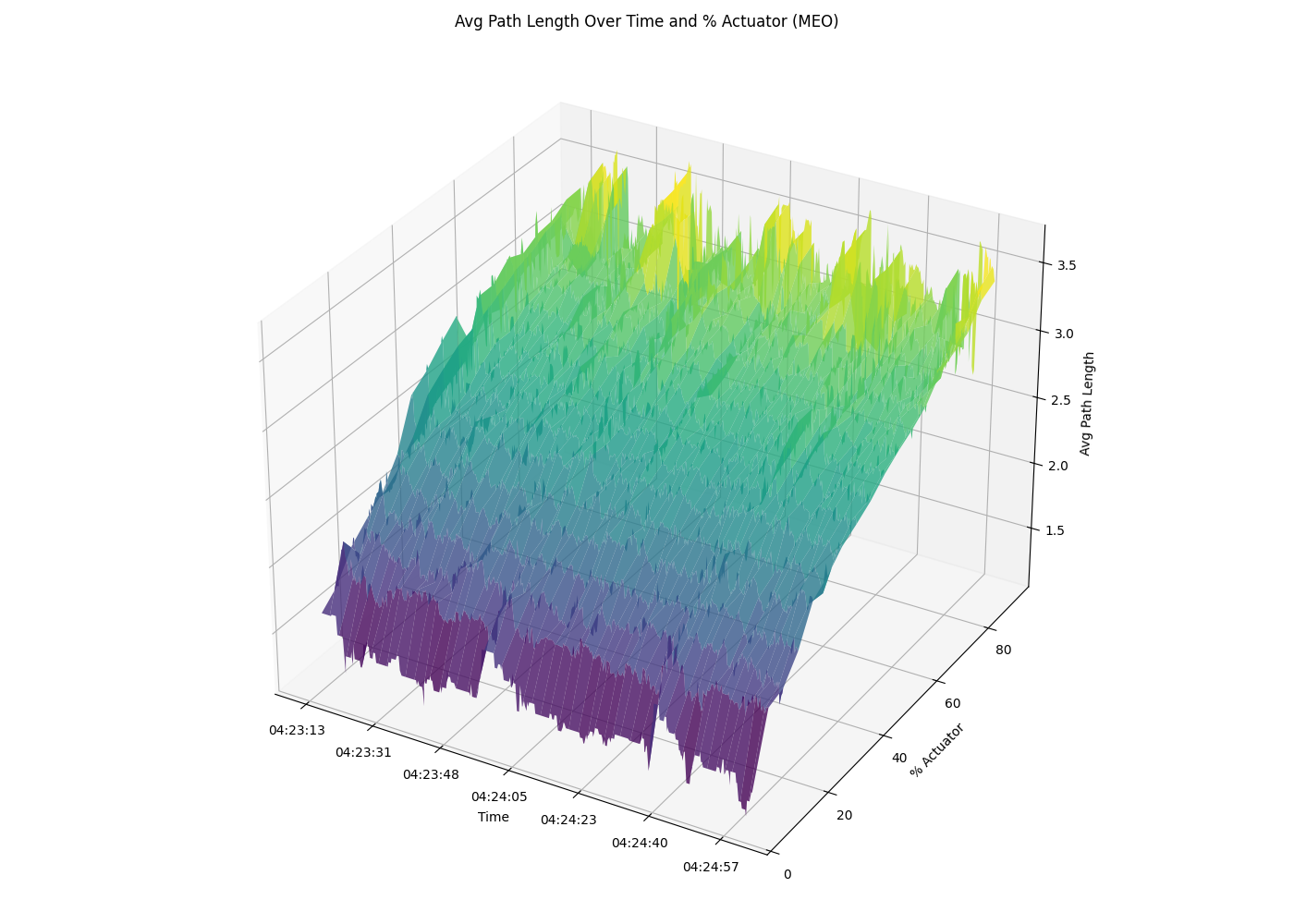}
\caption{}
\end{subfigure}
\caption{Performance metrics for MEO actuator network topology over time. (a) MST, (b) Degree variance, and (c) Average path length.}
\label{fig:meo2}
\end{figure}

\begin{figure}[h!]
\centering
\begin{subfigure}{0.31\linewidth}
  \centering
    \includegraphics[width=\linewidth, trim=0 0 0 40, clip]{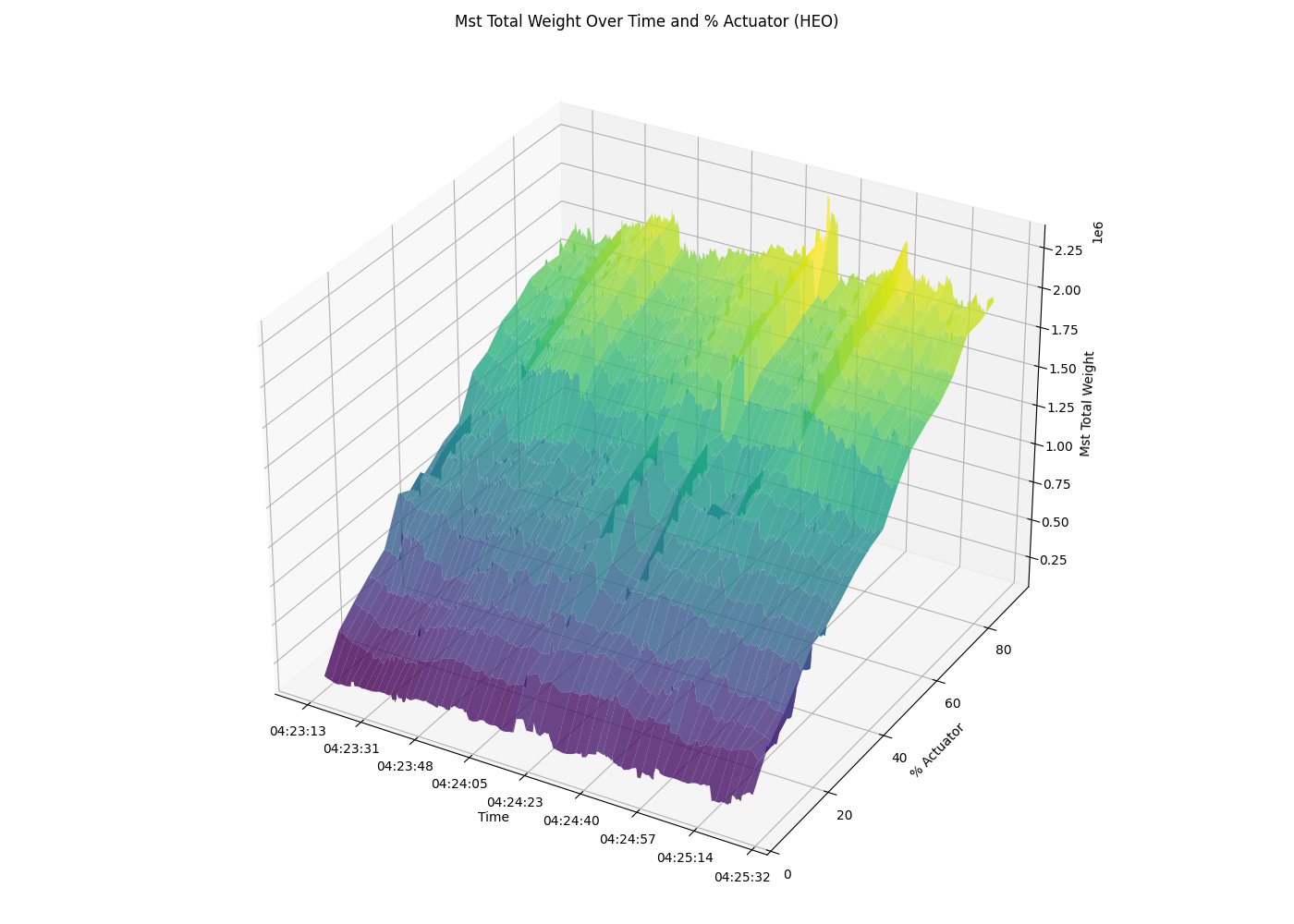}
\caption{}
\end{subfigure}%
\begin{subfigure}{0.31\linewidth}
  \centering
    \includegraphics[width=\linewidth, trim=0 0 0 40, clip]{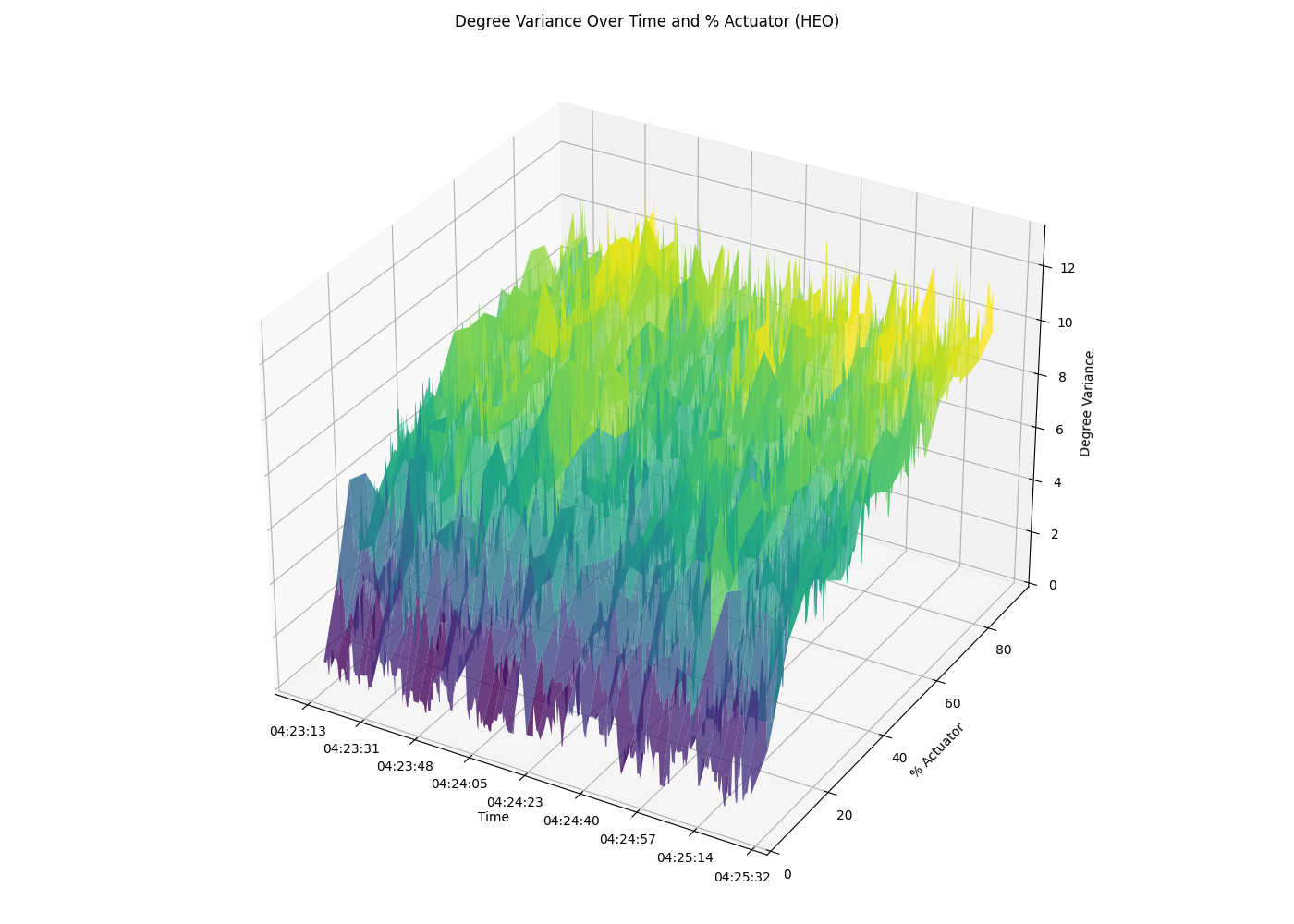}
\caption{}
\end{subfigure}%
\begin{subfigure}{0.31\linewidth}
  \centering
    \includegraphics[width=\linewidth, trim=0 0 0 40, clip]{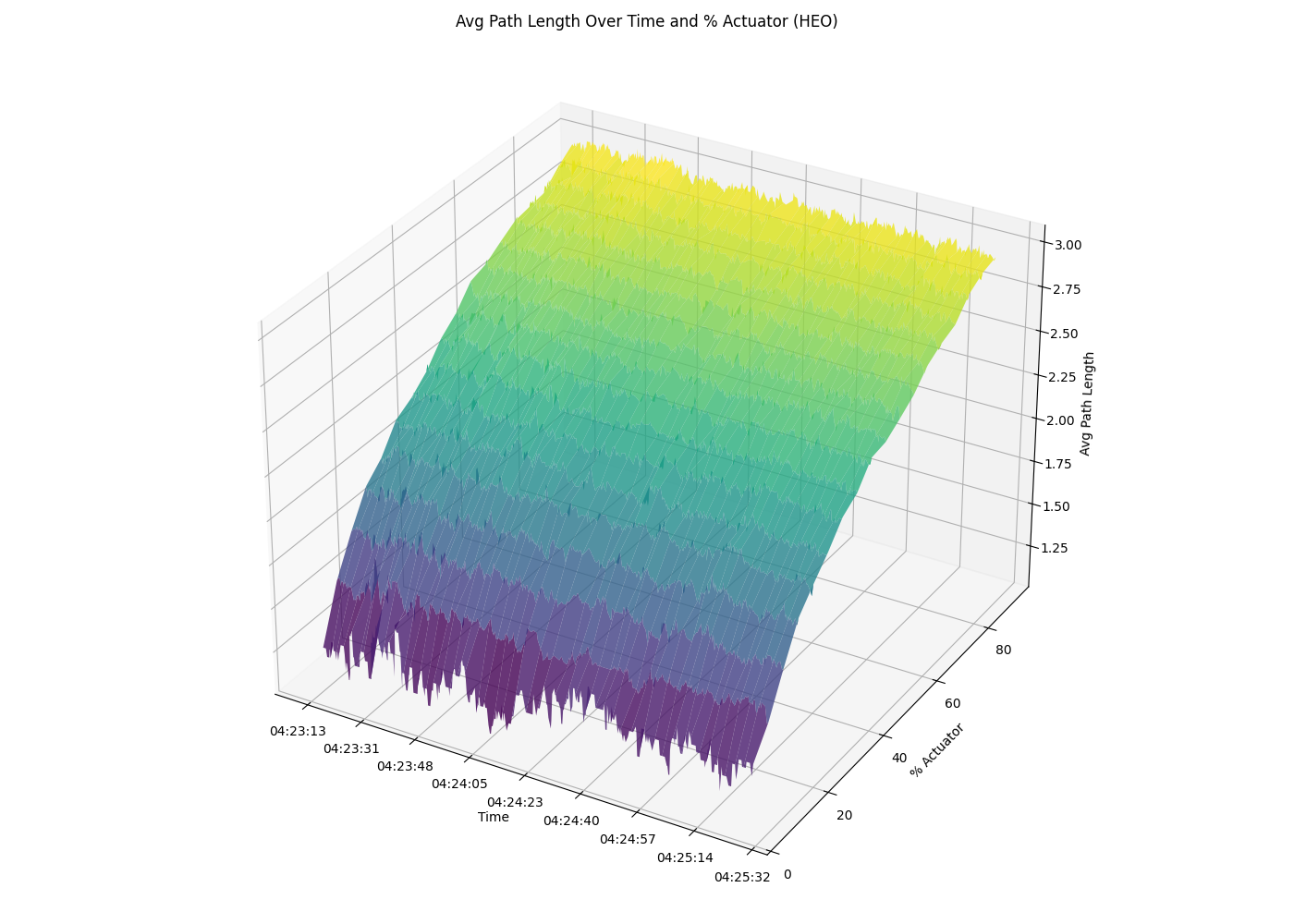}
\caption{}
\end{subfigure}%
\caption{Performance metrics for HEO actuator network topology over time. (a) MST, (b) Degree variance, and (c) Average path length.}
\label{fig:heo2}
\end{figure}

\begin{figure}[h!]
\centering
\begin{subfigure}{0.31\linewidth}
  \centering
    \includegraphics[width=\linewidth, trim=0 0 0 40, clip]{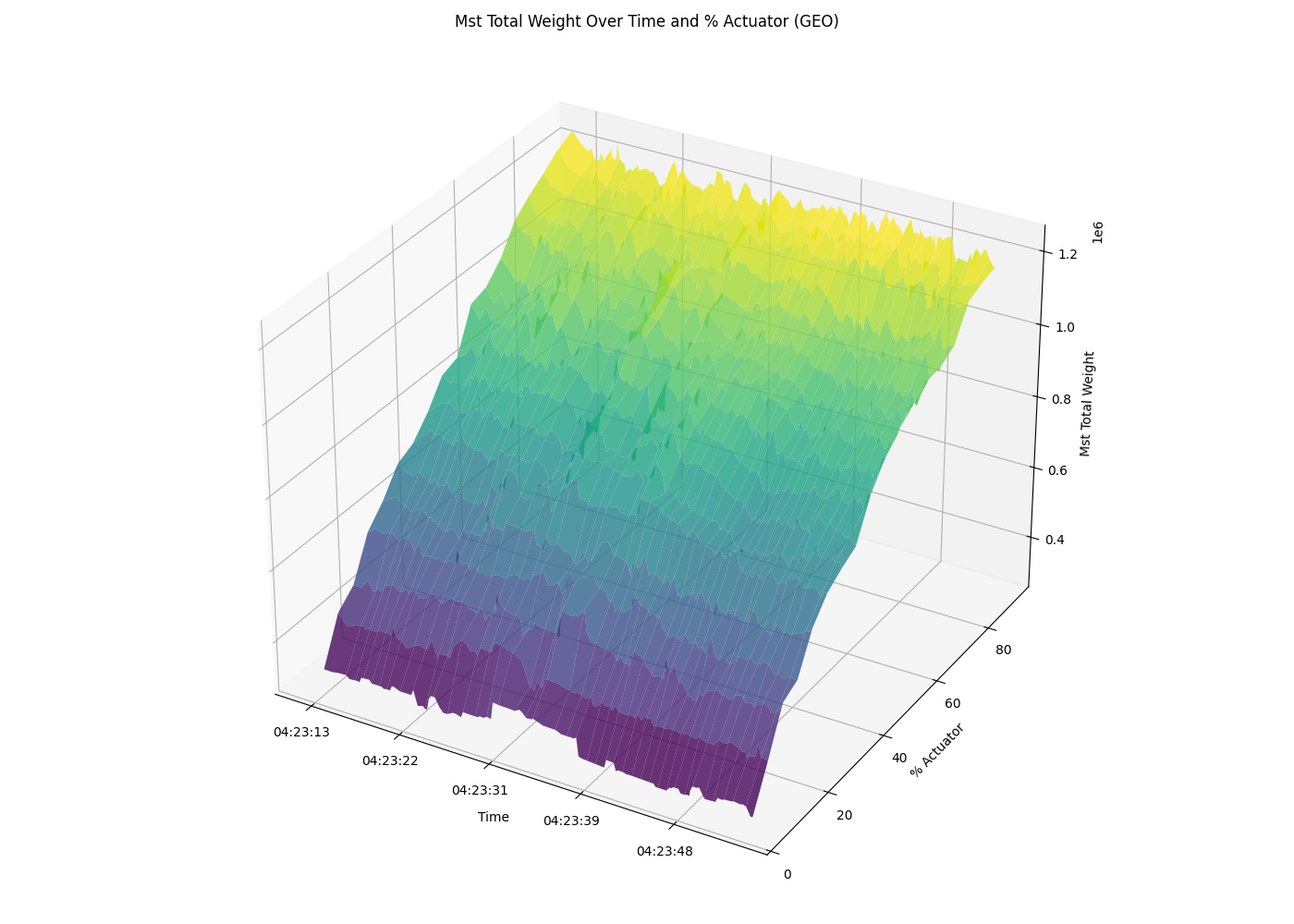}
\caption{}
\end{subfigure}%
\begin{subfigure}{0.31\linewidth}
  \centering
    \includegraphics[width=\linewidth, trim=0 0 0 40, clip]{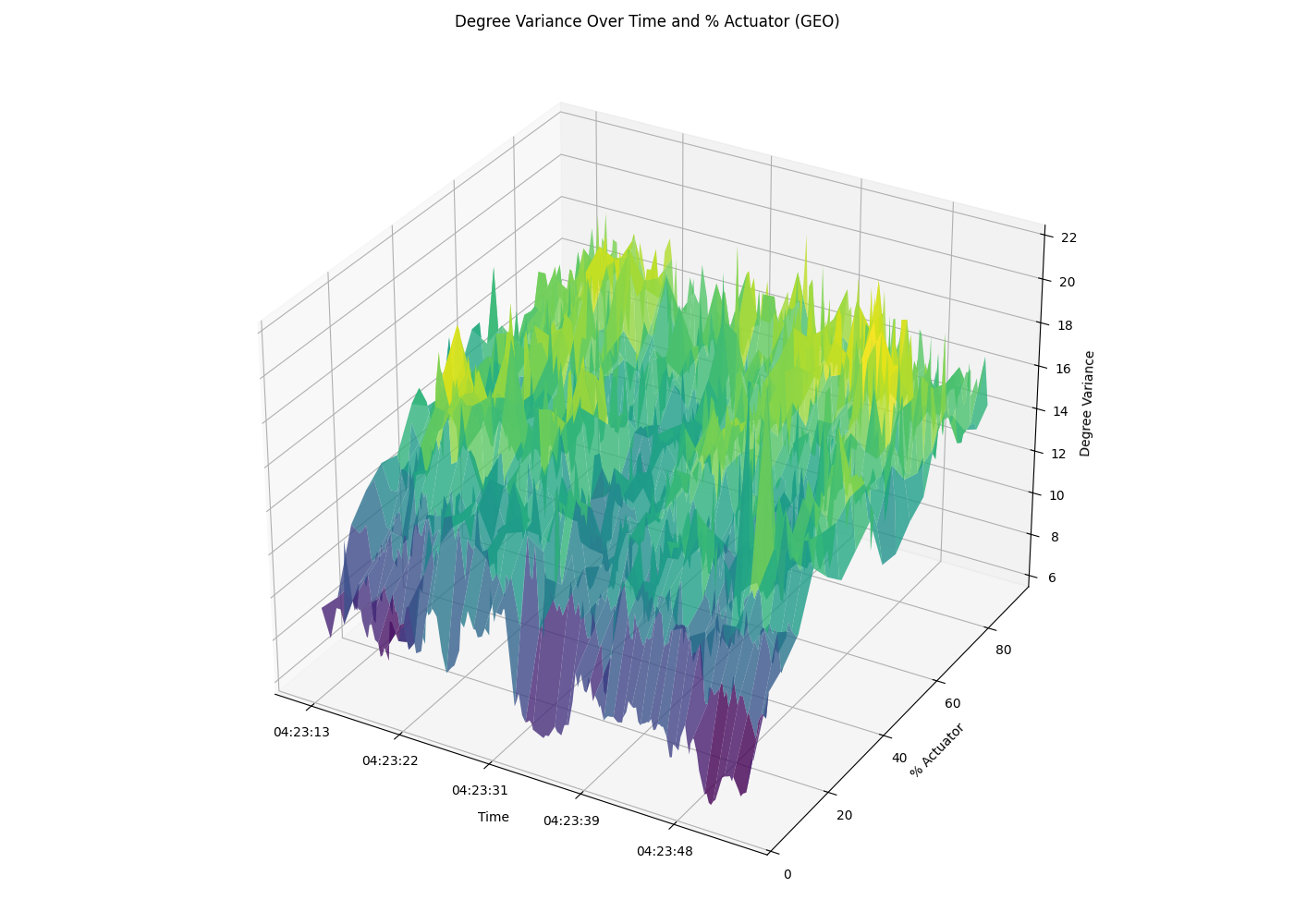}
\caption{}
\end{subfigure}%
\begin{subfigure}{0.31\linewidth}
  \centering
    \includegraphics[width=\linewidth, trim=0 0 0 40, clip]{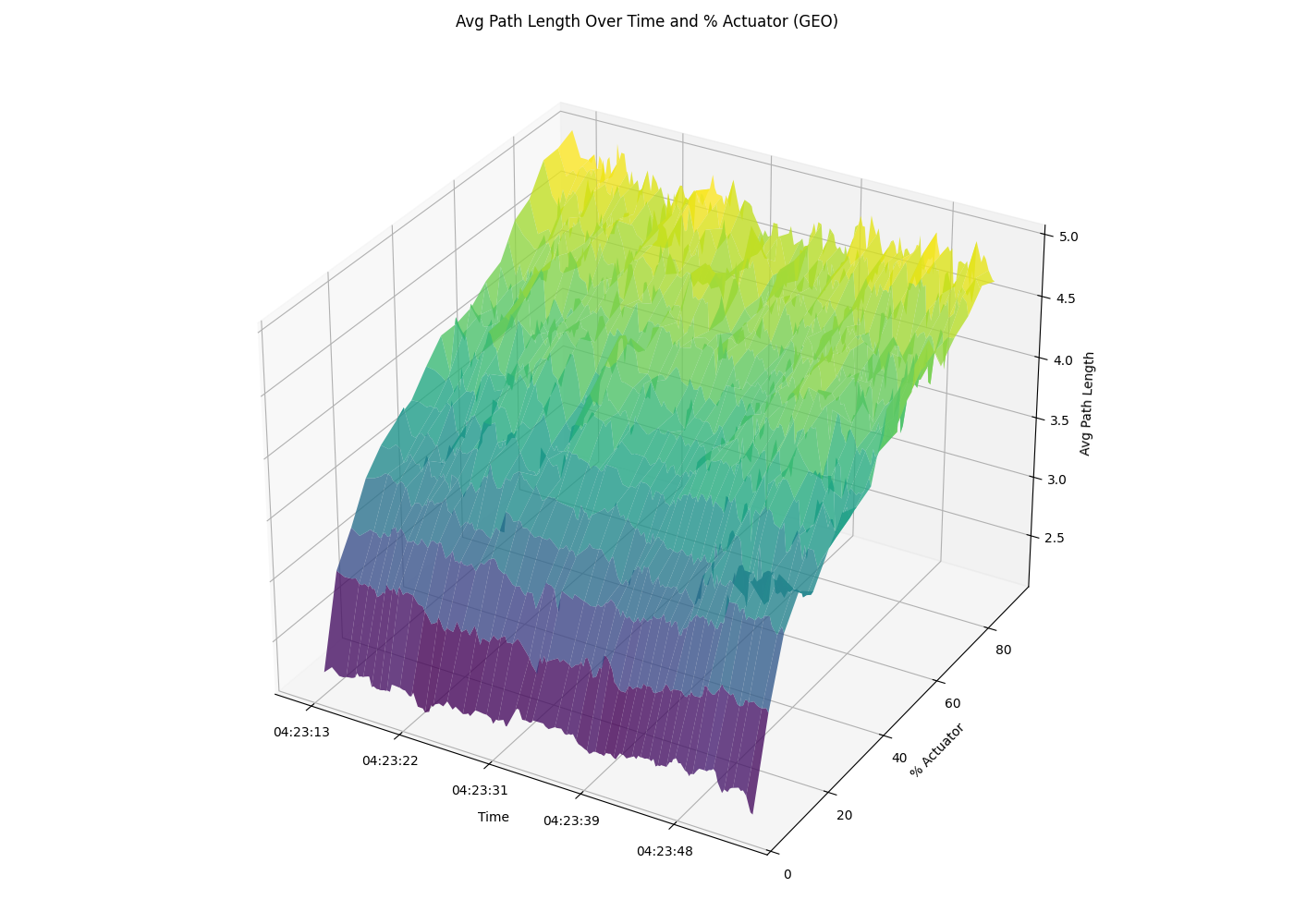}
\caption{}
\end{subfigure}
\caption{Performance metrics for GEO actuator network topology over time. (a) MST, (b) Degree variance, and (c) Average path length.}
\label{fig:geo2}
\end{figure}

Next, we report in Figs.~\ref{fig:leo2}, \ref{fig:meo2}, \ref{fig:heo2}, and \ref{fig:geo2} the performance metrics for the actuator network topology over time to show the MST, degree variance and average path length in 3D. For LEO, the MST shows significant variation but remains the lightest due to short distances between actuators, with fluctuations caused by bridging to farther nodes. Degree variance stabilizes despite occasional peaks, likely due to node velocity and their in-and-out movement in the frame and the average path length decreases as more actuators are added, benefiting from higher node density. In contrast, higher orbits (MEO, HEO, and GEO) show a consistent increase of trend in MST as actuator percentages increases, with the degree variance increasing, indicating less uniform connectivity, suggesting less efficient communication due to node dispersion. 

The results highlight a diminishing return in network efficiency with increasing actuators presence for higher orbits HEO, MEO and GEO. This is attributed to the elongated nature of ISLs, which lead to uneven spatial distribution and stretched communication paths as more nodes are introduced.  However, these orbits compensate for this with their larger coverage areas and slower-moving nodes, making them well-suited for managing broad regions of responsibility and can  still support reliable long-range communication, benefiting from the slower movement of nodes which results in more stable links (less variability in the degree). 

For synchronization potential, these higher orbits' actuators nodes, that are more \say{central}, help maintain stability in larger regions, providing consistent communication over longer periods, though at a slower pace. On the other hand, LEO offers a denser, more compact network, which allows for faster information propagation and greater scalability. This makes LEO ideal for tasks requiring rapid communication and responsiveness.

\section{Collective Behavior}

Based on the analysis of network behavior across various orbital regimes, we propose a hierarchical consensus framework tailored to the distinct characteristics and requirements of the regimes as the number of actuators increases. This framework extends our prior work in LEO \cite{benchoubane2025nextgenspacebasedsurveillanceblockchain}, where we highlighted the optimization of the consensus process through role segregation among actuators, aiming to enhance decision-making efficiency and reduce overall latency.

\subsection{Hierarchical Model: Approvers and Verifiers}
Algorithm~\ref{Hierarchical} presents the consensus process for the network where we define two main roles within the actuators network: approvers and verifiers.
\begin{enumerate}
    \item \textit{Approvers}: These actuators are responsible for approving data before it is passed to the verifiers for validation. Their main role is to ensure that the data is ready for verification, making them key in maintaining the flow of information in the network.
    \item \textit{Verifiers}: These actuators are responsible for verifying the data once it has been approved. Verifiers maintain the truthfulness of the network by ensuring that the approved data is accurate, consistent, and trustworthy.
\end{enumerate}

\begin{algorithm}[h!]
    \caption{Hierarchical consensus algorithm per phase.}
    \label{Hierarchical}

    \textbf{Input:} A set of actuators $A$, a set of verifiers $v$, and a set of approvers $a_p$.\\
    \textbf{Output:} Consensus result: "Verification is valid" or "Repeat Process!".

    \textbf{Step 1: Approval Phase (Approvers)}
    \begin{itemize}
        \item Each approver node selects data for approval from source nodes.
        \item Approver nodes check the validity of the data.
        \item Each approver node computes a signature $c_k$ for the data with its own key.
        \item Approvers release their computed signatures $c_k$ publicly for inclusion in the consensus process.
    \end{itemize}

    \textbf{Step 2: Verification Phase (Verifiers)}
    \begin{itemize}
        \item Each verifier node computes the sum of the approvers' signatures and confirms the result.
        \item Each verifier node then computes its own signature $d_k$ for the data.
        \item Verifiers release their computed signatures publicly.
    \end{itemize}

    \textbf{Step 3: Consensus Check}
    \begin{itemize}
        \item If the sum of the approvers' signatures matches the sum of the verifiers' signatures:
        \item If they match, output: "Verification is valid!"
        \item If they do not match, output: "Repeat Process!"
    \end{itemize}
\end{algorithm}

To guide the distribution and assignment of roles, we differentiate between lower and higher orbital regimes.

\subsection{Lower Orbital Regimes}

\begin{figure}[h!]
    \centering
    \includegraphics[width=0.8\linewidth]{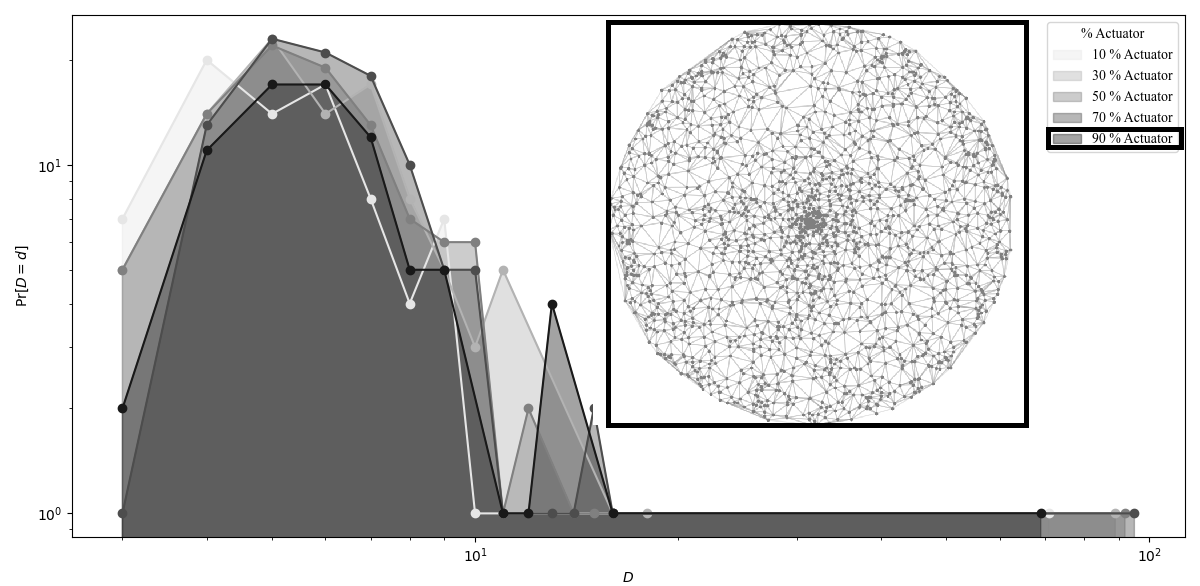}
    \caption{Degree distribution in a LEO  network with varying actuator percentages.}
    \label{fig:leocon}
\end{figure}

In LEO networks, performance improves significantly as the number of actuators increases and the network exhibits a degree distribution that approximates a Gaussian shape, as shown in Fig.~\ref{fig:leocon}. Nodes with higher degree, which lie at the peak of this distribution, have a higher number of connections and are thus well-suited to handle both the roles of approvers and verifiers. These located nodes, due to their broad connectivity, are ideal for efficiently managing the data flow and ensuring that data is both approved and verified. However, instead of splitting and assigning roles to specific nodes permanently, the role assignment should be based on the node's current connectivity, which may change as nodes move in their orbits. This way, we can also accommodate the shifting nature of LEO, where responsibility regions evolve over time.

\subsection{Higher Orbital Regimes}

\begin{figure}[h!]
    \centering
    \includegraphics[width=0.3\linewidth]{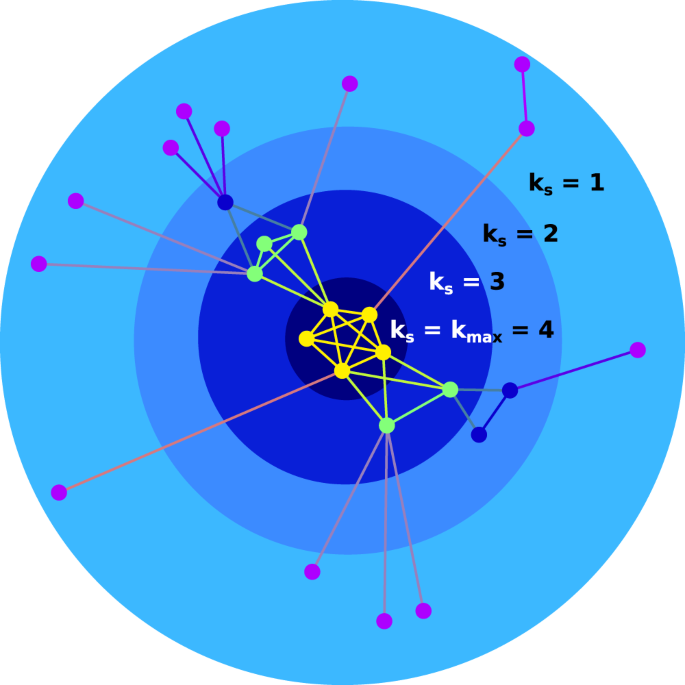}
    \caption{Sketch of $k$-core decomposition \cite{Burleson-Lesser2020}.}
    \label{fig:kcoredecomposition}
\end{figure}

In higher orbital regimes,  the challenges are different due to the ISLs and the stretched communication paths as more nodes are added. The spatial distribution becomes uneven, and thus, the approach to role assignment needs adjustment. Here, the central nodes should still dictate the role distribution, but this time, we can rely on $k$-core decomposition to identify the central actuators. This decomposition, as shown in Fig.~\ref{fig:kcoredecomposition}, identifies the core nodes responsible for the network’s connectivity. In these regimes, the verifiers would be the central core nodes, while approvers would occupy nodes in the next lower tier, those identified by the $k$ and $k-1$ level of the k-core decomposition retrospectively. 

\noindent While the hierarchical consensus framework provides a strong foundation, future work will focus on validating these role assignment algorithms. We will investigate the optimal percentage of nodes in each role, the dynamics of role assignment as the network evolves, and the impact of the k-core value on network robustness. Moreover, the effectiveness of the consensus process will be evaluated using service metrics that consider actuator performance, link weights, and network redundancy.





\section{Conclusions} 
In conclusion, this paper examines the potential for distributed SDA through on-orbit collaboration among heterogeneous satellites. We introduce a graph-theoretic approach using Voronoi tessellations to define spatial responsibilities and Delaunay triangulations to model communication paths between actuator nodes. Our analysis, applied to LEO, MEO, HEO, and GEO, demonstrates the feasibility of leveraging on-orbit assets distributed SDA, addressing the challenges posed by the increasing number of space assets and the need for dynamic, decentralized decision-making. 

\begin{acknowledgment}
	This work was supported in part by the Tier 1 Canada Research Chair program. The authors would also like to express sincere gratitude to the Canadian Space Agency for their support in funding participation for the 18th International Conference on Space Operations under the Announcement of Opportunity (AO) of Canadian Student Participation in Space Conferences and Training Events 2025.
\end{acknowledgment}



\bibliographystyle{IEEEtran}
\bibliography{references}

\end{document}